\documentclass[letter,twoside,english,showkeys]{revtex4}
\usepackage[T1]{fontenc}
\usepackage[latin1]{inputenc}
\usepackage{amsmath}
\usepackage{graphicx}
\usepackage{amssymb}

\makeatletter

\providecommand{\tabularnewline}{\\}

\newcommand{\eg}{e.g.\ }
\newcommand{\ie}{i.e.\ }

\newcommand{\sect}[1]{Sec.~\ref{Sec:#1}}      \newcommand{\Lsect}[1]{\label{Sec:#1}}

\newcommand{\eps}{\varepsilon}
\newtheorem{thm}{Theorem}[section]

\newtheorem{mytheorem}[thm]{Theorem}
\newcommand{\field}[1]{\mathbb{#1}}

\newcommand{\NN}{\field{N}}

\newcommand{\Pre}{\mathrm{Pre}}
\newcommand{\Post}{\mathrm{Post}}

\newcommand{\rank}{\mathrm{rank}}
\newcommand{\tauh}{\frac{\tau}{2}}

\newcommand{\mysubstack}[2]{\substack{#1\\#2} }
\newcommand{\reset}{\rightarrow}

\newcommand{\hD}[2]{U^{-1}(U(#1+D_{i,#2})+
  \eps_{i j_{#2}})=:\beta_{i,#2}}

\usepackage{babel}
\makeatother
\begin{document}

\title{The simplest problem in the collective dynamics of neural networks:\\
 Is synchrony stable?}

\author{Marc Timme and Fred Wolf}

\affiliation{Max Planck Institute for Dynamics and Self-Organization,\\
Bernstein Center for Computational Neuroscience, and\\
University of Göttingen, Bunsenstr. 10, 37073 Göttingen, Germany}

\begin{abstract}
For spiking neural networks we consider the stability problem of global
synchrony, arguably the simplest non-trivial collective dynamics in
such networks. We find that even this simplest dynamical problem --
local stability of synchrony -- is non-trivial to solve and requires
novel methods for its solution. In particular, the discrete mode of
pulsed communication together with the complicated connectivity of
neural interaction networks requires a non-standard approach. The
dynamics in the vicinity of the synchronous state is determined by
a multitude of linear operators, in contrast to a single stability
matrix in conventional linear stability theory. This unusual property
qualitatively depends on network topology and may be neglected for
globally coupled homogeneous networks. For generic networks, however,
the number of operators increases exponentially with the size of the
network. 

We present methods to treat this multi-operator problem exactly.
First, based on the Gershgorin and Perron-Frobenius theorems, we derive
bounds on the eigenvalues that provide important information about
the synchronization process but are not sufficient to establish the
asymptotic stability or instability of the synchronous state.  We
then present a complete analysis of asymptotic stability for topologically
strongly connected networks using simple graph-theoretical considerations.

For inhibitory interactions between dissipative (leaky) oscillatory
neurons the synchronous state is stable, independent of the parameters
and the network connectivity. These results indicate that pulse-like interactions play a profound role in network dynamical systems, and in particular in the dynamics of biological synchronization, unless the coupling is homogeneous and all-to-all. The concepts introduced here are expected to also facilitate
the exact analysis of more complicated dynamical network states, for instance the irregular balanced activity in cortical neural networks.
\end{abstract}


\maketitle

\section{Introduction}

\label{Sec:Intro}

Understanding the collective dynamics of biological neural networks
represents one of the most challenging problems in current theoretical
biology. For two distinct reasons, the study of synchronization in
large neuronal networks plays a paradigmatic role in theoretical studies
of neuronal dynamics. Firstly, synchrony is an ubiquitous collective
feature of neural activity. Large-scale synchronous activity has been
observed by spatially coarse methods such as electro- or magnetoencephalography.
Complementing experiments of parallel recordings of spiking activity
of individual cells have shown that the synchronization of firing
times of indiviual units is often precise with a temporal scatter
of the order of a few milliseconds \cite{Eckhorn:1988:121,Gray:1989:334}.
This precise locking has been observed over significant distances
in the cortex and even across hemispheres \cite{Keil:1999:7152}.
Synchronous activity also plays an important role in pathological
state such as epileptic seizures \cite{Gray:1994:11}. Secondly, the
synchronous state is arguably the simplest non-trivial collectively
coordinated state of a network dynamical system. Mathematically, it
is therefore one of the most throroughly investigated states in the
dynamics of biological neural networks. Following the paradigm set
by the seminal works of Winfree and Kuramoto on the dynamics of biological
oscillators, most studies of synchronization have utilized either
temporal or population averaging techniques to map the pulse-coupled
dynamics of biological neural networks to effective models of phase-coupled
oscillators or density dynamics \cite{Abbott:1993:1483,Herz:1995:1222,Tsodyks:1993:1280,Corall:1995:3697,Gerstner:1996:1653,vanVreeswijk:1996:5522,vanVreeswijk:2000:5110,Bressloff:1997:2791,Bressloff:1997:4665,Golomb:2000:1095,Hansel:2001:4175,Brunel:1999:1621,Ernst:1995:1570,Ernst:1998:2150,Mirollo:1990:1645,Peskin:1984,Senn:2001:1143,Roxin:2005:238103}. 

While this approach has proven to be very informative, it provides
exact results only under highly restrictive conditions. In fact, studies
approaching the dynamics of biological neural networks using exact
methods \cite{Tsodyks:1993:1280,Timme:2002:154105,Timme:2003:1,Zumdieck:2002,Zillmer:2006:036203,Zillmer:2007:1960,Zillmer:2007:046102,Memmesheimer:2006:PD06,Memmesheimer:2006:188101,Timme:2004:074101,Timme:2006:015108,Wu:2007:789,Denker:2002:074103},
have over the past decade revealed numerous examples of collective
behaviours that obviouslly are outside the standard repertoir of behaviours
expected from a smoothly coupled dynamical system. These include several
new and unexpected phenomena such as unstable attractors
\cite{Timme:2002:154105,Timme:2003:1,Ashwin:2005:2035}, stable chaos
\cite{Timme:2002:154105,Timme:2003:1,Jahnke:2007:1}, topological
speed limits to coordinating spike times \cite{Timme:2004:074101,Timme:2006:015108},
and extreme sensitivity to network topology \cite{Zumdieck:2004:244103}.
The occurence of these phenomena may signal that the proper theoretical
analysis of collective neuronal dynamics mathematically represents
a much more challenging task than is currently appreciated. 

In this article, we study the impact of pulse-coupling, delayed interactions
and complicated network connectivity on the exact microscopic dynamics
of neural networks and expose the mathematical complexity that emerges
already when considering the seemingly simple problem of neuronal
synchronization. Utilizing the Mirollo-Strogatz phase representation
of individual units, we present an analytical treatment of finite
networks of arbitrary connectivity. The results obtained unearth an
unanticipated subtlety of the nature of this stability problem: It
turns out that a single linear operator is not sufficient to represent
the local dynamics in these systems. Instead, a large number of linear
operators, depending on rank order of the perturbation vector, is
needed to represent the dynamics of small perturbations. We present
methods to characterize the eigenvalues of all operators arising for
a given network. Universal properties of the stability operators lead
to exact bounds on the eigenvalues of all operators and also provide
a simple way to prove plain but not asymptotic stability for any network.
We then show for topologically strongly connected networks that asymptotic
stability of the synchronized state can be demonstrated by graph-theoretical
considerations. We find that for inhibitory interactions the synchronous
state is stable, independent of the parameters and the network connectivity.
A part of this work that considers plain (non-asymptotic) stability
without networks constraints has been briefly reported before for
the case of inhibitory coupling \cite{Timme:2002:258701}.These results 
indicate that pulse-like interactions play a profound role in network 
dynamical systems, and in particular in the dynamics of biological 
synchronization, unless the coupling is homogeneous and all-to-all.
They highlight the need for exact mathematical tools that can handle the dynamic complexity of neuronal systems at the microscopic level.

\section{The Phase Representation of Pulse-Coupled Networks}

\label{Sec:MSmodel}

We consider networks of $N$ pulse-coupled oscillatory units, neurons, with
delayed interactions. A phase-like variable $\phi_{j}(t)\in(-\infty,1]$
specifies the state of each oscillator $j$ at time $t$ such that
the difference between the phases of two oscillators quantifies their
degree of synchrony, with identical phases for completely synchronous
oscillators. The free dynamics of oscillator $j$ is given by \begin{equation}
\frac{d\phi_{j}}{dt}=1.\label{eq:phi_dot_1}\end{equation}
Whenever oscillator $j$ reaches a threshold \begin{equation}
\phi_{j}(t)=1\label{eq:threshold_phi}\end{equation}
the phase is reset to zero \begin{equation}
\phi_{j}(t^{+})=0\label{eq:reset_phi}\end{equation}
and a pulse is sent to all other post-synaptic oscillators $i\in\Post(j)$.
These oscillators $i$ receive this signal after a delay time $\tau$.
The interactions are mediated by a function $U(\phi)$ specifying
a 'potential' of an oscillator at phase $\phi.$ The function $U$
is assumed twice continuously differentiable, monotonically increasing,
\begin{equation}
U'>0,\end{equation}
concave (down), \begin{equation}
U''<0,\end{equation}
and normalized such that \begin{equation}
U(0)=0\quad\textrm{ and }\quad U(1)=1.\end{equation}

\begin{figure}[!tph]
\begin{centering}\includegraphics[width=6cm]{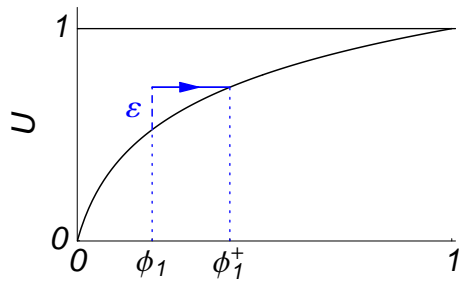}\par\end{centering}

\caption{An incoming pulse of strength $\eps$ induces a phase jump $\phi_{1}^{+}:=\phi_{1}(t^{+})=U^{-1}(U(\phi_{1}(t))+\eps)=H_{\eps}(\phi_{1})$
that depends on the state $\phi_{1}:=\phi_{1}(t)$ of the oscillator
at time $t$ of pulse reception. Due to the curvature of $U$, an
excitatory pulse ($\eps>0$) induces an advancing phase jump. If the
incoming pulse puts the potential above threshold ($U(\phi_{1})+\eps>1$),
the phase is reset to zero ($\phi_{1}^{+}=0$). An inhibitory pulse
($\eps<0$) would induce a regressing phase jump such that the phase
may assume negative values (not shown).\label{Fig:potential_to_phase} }
\end{figure}

For a general $U(\phi)$ we define the transfer function \begin{equation}
H_{\eps}(\phi)=U^{-1}(U(\phi)+\eps)\label{eq:H}\end{equation}
that represents the response of an oscillator at phase $\phi$ to
an incoming sub-threshold pulse of strength $\eps$. Depending on
whether the input $\eps_{ij}$ received by oscillator $i$ from $j$
is sub-threshold, \begin{equation}
U(\phi)+\eps_{ij}<1,\end{equation}
or supra-threshold, \begin{equation}
U(\phi)+\eps_{ij}\geq1,\end{equation}
the pulse sent at time $t$ (Eq.~(\ref{eq:threshold_phi})) induces
a phase jump after a delay time $\tau$ at time $t+\tau$ according
to \begin{equation}
\phi_{i}((t+\tau)^{+})=\left\{ \begin{array}{lll}
H_{\eps_{ij}}(\phi_{i}(t+\tau)) & \mathrm{if} & U(\phi_{i}(t+\tau))+\eps_{ij}<1\\
0 & \mathrm{if} & U(\phi_{i}(t+\tau))+\eps_{ij}\geq1\end{array}\right..\label{eq:phase_jump}\end{equation}
In the second case also a pulse is emitted by oscillator $i$ . The
phase jump (Fig.~\ref{Fig:potential_to_phase}) depends on the phase
$\phi_{i}(t+\tau)$ of the receiving oscillator $i$ at time $t+\tau$
after the signal by oscillator $j$ has been sent at time $t$, the
effective coupling strength $\eps_{ij}$, and the nonlinear potential
function $U$. %
 The interaction from unit $j$ to unit $i$ is either excitatory
($\varepsilon_{ij}>0$) inducing advancing phase jumps (cf.~Fig.~\ref{Fig:potential_to_phase})
or inhibitory ($\varepsilon_{ij}<0$) inducing retarding phase jumps.
If there is no interaction from $j$ to $i$, we have $\eps_{ij}=0$
and . There are two immediate differences between inhibitory and excitatory
inputs. First, while inhibitory input $\eps_{ij}<0$ is always sub-threshold,
excitatory input $\eps_{ij}>0$ may also be supra-threshold and thus
induce an instantaneous reset to zero according to Eq.~(\ref{eq:phase_jump}).
Second, in response to the reception of an inhibitory pulse, the phases
of the oscillators may assume negative values whereas for excitatory
coupling they are confined to the interval $[0,1]$. We remark here
that we do not consider the dynamics (and perturbations) of variables
that encode the spikes in transit, i.e. spikes that have been sent
but not yet received at any instant of time. These additional variables
make the system formally higher-dimensional. However, under the conditions
considered below, in particular for networks of identical neurons
with inhibitory interactions, earlier rigorous work \cite{Ashwin:2005:2035}
shows that these spike time variables lock to the phase variables
once all spikes present in the system initially have arrived and every
neuron has emitted at least one spike thereafter (which takes finite
total time). Thus, for our purposes, we consider the dynamics and
perturbations of the phase variables only.

Choosing appropriate functions $U$ the dynamics of a wide variety
of pulse coupled systems can be represented. In particular, any differential
equation for an individual neural oscillator of the form \begin{equation}
{\dot{V}_{i}}=f(V_{i})+S_{i}(t)\label{eq:Vdot}\end{equation}
together with the threshold firing condition \begin{equation}
V_{i}(t)=1\,\Rightarrow V_{i}(t^{+})=0\label{eq:Vthreshold}\end{equation}
where \begin{equation}
S_{i}(t)=\sum_{j,m}\eps_{ij}\delta(t-t_{j.m})\label{eq:networkcurrentS}\end{equation}
is the synaptic current from the network and $t_{j,m}$ is the time
of the $m^{th}$ firing of neuron $j$. As long as the free ($S_{i}(t)\equiv0$)
dynamics has a periodic solution $V(t)$ with period $T$ and negative
curvature, the function $U$ can be taken as the scaled solution,
\begin{equation}
U(\phi):=V(\phi T).\label{eq:UV}\end{equation}
Thus a general class of pulse-coupled oscillator networks, defined
by Eq.~(\ref{eq:Vdot}-\ref{eq:networkcurrentS}), can be mapped
onto the normalized phase representation (Eqs.~(\ref{eq:phi_dot_1})--(\ref{eq:phase_jump})).

In this paper, we consider a system of $N$ coupled Mirollo-Strogatz
oscillators which interact on directed graphs by sending and receiving
pulses. The structure of this graph is specified by the sets $\Pre(i)$
of presynaptic oscillators that send pulses to oscillator $i$. For
simplicity we assume no self-interactions, $i\notin\Pre(i)$. The
sets $\Pre(i)$ completely determine the topology of the network,
including the sets $\Post(i)$ of postsynaptic oscillators that receive
pulses from $i$. The coupling strength between oscillator $j$ and
oscillator $i$ is given by $\eps_{ij}$ such that\begin{equation}
\begin{array}{ll}
\eps_{ij}\neq0 & \textrm{if }j\in\Pre(i)\\
\eps_{ij}=0 & \textrm{otherwise}.\end{array}\end{equation}
Thus a connection is considered to be present if the connection strength
is non-zero. All analytical results presented are derived for the
general class of interaction functions $U$ introduced above. The
structure of the network is completely arbitrary except for the restriction
that every oscillator has at least one presynaptic oscillator.

\section{regular and irregular dynamics in pulse-coupled networks}

\label{sec:Coexistence}

\begin{figure}[!tph]
\begin{centering}\includegraphics[width=10cm]{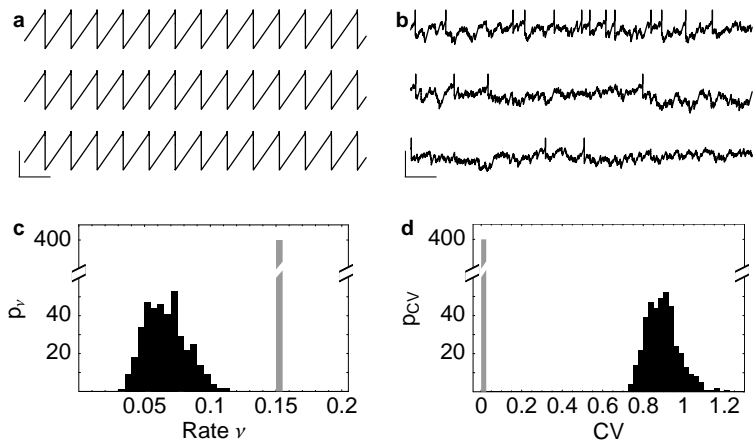}\par\end{centering}

\caption{Coexistence of (a) synchronous and (b) irregular dynamics in a random
network ($N=400$, $p=0.2$, $I=4.0,$ $\eps=-16.0$, $\tau=0.035$),
cf.~\cite{Timme:2002:258701}. (a),(b): Trajectories of the potential
$U(\phi_{i})$ of three oscillators (angular bars: time scale (horizontal)
$\Delta t=8$; potential scales (vertical) (a) $\Delta U=8$, (b)
$\Delta U=2$ ; spikes of height $\Delta U=1$ added at firing times).
(c),(d): Distributions (c) $p_{\nu}$ of rates and (d) $p_{\mathrm{CV}}$
of the coefficient of variation, displayed for the irregular (dark
gray) and synchronous (light gray) dynamics. Figure modified from
\cite{Timme:2002:258701}. \label{fig:trajectoryISIcv}}
\end{figure}

\begin{figure}[!tph]
\begin{centering}\includegraphics[width=9cm]{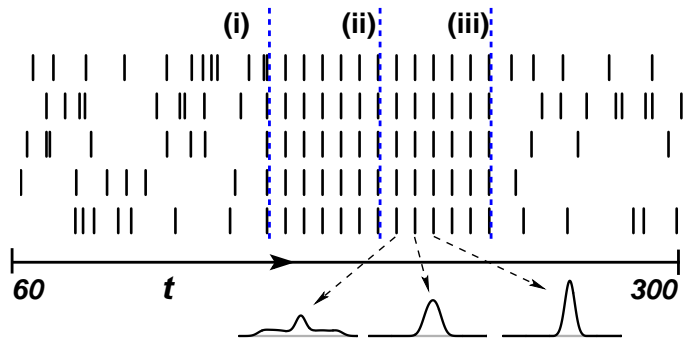}\par\end{centering}

\caption{Switching between irregular and synchronous dynamics ($N=400$, $p=0.2$,
$I=4.0$, $\eps=-16.0$, $\tau=0.14$). Firing times of five oscillators
are shown in a time window $\Delta t=240$. Vertical dashed lines
mark external perturbations: (i) large excitatory pulses lead to synchronous
state, (ii) a small random perturbation ($|\Delta\phi_{i}|\leq0.18$)
is restored, (iii) a sufficiently large random perturbation ($|\Delta\phi_{i}|\leq0.36$)
leads to an irregular state. Bottom: Time evolution of the spread
of the spike times after perturbation (ii), total length $\Delta t=0.25$
each. Figure modified from \cite{Timme:2002:258701}. \label{fig:spike_perturb}}
\end{figure}

The synchronous state in which

\begin{equation}
\phi_{i}(t)=\phi_{0}(t)\quad\textrm{for all }i\end{equation}
is arguably one of the simplest ordered states a network of pulse-coupled
oscillators may assume. This synchronous state exists if and only
if the coupling strengths are normalized such that \begin{equation}
\sum_{j\in\Pre(i)}\eps_{ij}=\eps.\label{eq:epsnormalized}\end{equation}
One should keep in mind that this state whether stable or not is typically
not a global attractor in complex networks. To illustrate this point
let us briefly consider a specific example. Figures (Fig.\ \ref{fig:trajectoryISIcv}
and \ \ref{fig:spike_perturb}) show numerical results from a in
a randomly connected network of integrate-and-fire oscillators where
$U(\phi)=U_{\textrm{IF}}(\phi)=I(1-e^{-\phi T_{\textrm{IF}}})$ and
$T_{\textrm{IF}}=\ln(I/(I-1))$ with strong interactions. This network
exhibits a balanced state (cf.\ \cite{Brunel:1999:1621,vanVreeswijk:1996:1724,vanVreeswijk:1998:1321})
characterized by irregular dynamics. In this balanced state, found
originally in binary neural networks \cite{vanVreeswijk:1996:1724,vanVreeswijk:1998:1321},
inhibitory and excitatory inputs cancel each other on average but
fluctuations lead to a variability of the membrane potential and a
high irregularity in firing times (see also \cite{Brunel:1999:1621}).
Figures \ref{fig:trajectoryISIcv}a,b display sample trajectories
of the potentials $U(\phi_{i})$ of three oscillators for the same
random network, making obvious the two distinct kinds of coexisting
dynamics. Whereas in the synchronous state all oscillators display
identical periodic dynamics, in the balanced state oscillators fire
irregularly, asynchronous, and in addition differ in their firing
rates.%

The latter dynamical difference is quantified by a histogram $p_{\nu}$
of oscillator rates (Fig.\ \ref{fig:trajectoryISIcv}c) \begin{equation}
\nu_{i}=\left(\left\langle t_{i,n+1}-t_{i,n}\right\rangle _{n}\right)^{-1},\label{eq:ISIdefinition}\end{equation}
 the reciprocal values of the time averaged inter-spike intervals.
Here the $t_{i,n}$ are the times when oscillator $i$ fires the $n^{\mathrm{th}}$
time. The temporal irregularity of the firing-sequence of single oscillators
$i$ is measured by the coefficient of variation \begin{equation}
CV_{i}=\left(\nu_{i}^{2}\left\langle (t_{i,n+1}-t_{i,n})^{2}\right\rangle _{n}-1\right)^{\frac{1}{2}},\label{eq:CVdefinition}\end{equation}
 defined as the ratio of the standard deviation of the inter-spike
intervals and their average. A histogram $p_{\mathrm{CV}}$ of the
$CV_{i}$ shows that the irregular state exhibits coefficients of
variation near one, the coefficient of variation of a Poisson process.
Such irregular states occur robustly when changing parameters and
network topology. 

Figure \ref{fig:spike_perturb} illustrates the bistability of these
qualitatively different states by switching the network dynamics form
one state to the other by external perturbations. A simple mechanism
to synchronize oscillators that are in a state of irregular firing
is the delivery of two sufficiently strong external excitatory (phase-advancing)
\textit{pulses} that are separated by a time $\Delta t\in(\tau,1)$,
cf.\ Fig.\ \ref{fig:spike_perturb}. The first pulse then leads
to a synchronization of phases due to simultaneous supra-threshold
input, cf. (\ref{eq:H}). If there are traveling signals that have
been sent but not received at the time of the first pulse, a second
pulse after a time $\Delta t>\tau$ is needed that synchronizes the
phases after all internal signals have been received. In this network
the synchronous state is not affected by \textit{small random perturbations},
whereas \textit{large random perturbations} lead back to irregular
dynamics (Fig.\ \ref{fig:spike_perturb}). These features clearly
suggests that the synchronized states is an asymptotically stable
local attractor. However, with the exception of special cases no exact
treatment of this proposition exists. Below we will expose that a
general treatment of the stability of apparently simple synchronous
state reveals an unexpectedly complex mathematical setting.

\section{Constructing Stroboscopic Maps\label{sec:MapConstruction}}

\label{sec:strobmap}

We perform a stability analysis of the synchronous state \begin{equation}
\phi_{i}(t)=\phi_{0}(t)\quad\textrm{for all }i\end{equation}
 in which all oscillators display identical phase-dynamics $\phi_{0}(t)$
on a periodic orbit such that $\phi_{0}(t+T)=\phi_{0}(t)$. The period
of the synchronous state is given by \begin{equation}
T=\tau+1-\alpha\label{eq:period}\end{equation}
 where \begin{equation}
\alpha=U^{-1}(U(\tau)+\eps).\label{eq:alpha}\end{equation}
For simplicity, we consider the cases where the total input $\eps$
is sub-threshold, $U(\tau)+\eps<1$ such that $\alpha<1$. A perturbation
\begin{equation}
\boldsymbol{\delta}(0)=:\boldsymbol{\delta}=(\delta_{1},\ldots,\delta_{N})\end{equation}
 to the phases is defined by \begin{equation}
\delta_{i}=\phi_{i}(0)-\phi_{0}(0)\,.\label{eq:define_delta}\end{equation}
 If we assume that the perturbation is small in the sense that \begin{equation}
\max_{i}\delta_{i}-\min_{i}\delta_{i}<\tau\label{eq:taubounds}\end{equation}
this perturbation can be considered to affect the phases of the oscillators
at some time just after all signals have been received, i.e.~after
a time $t>t_{0}+\tau$ if all oscillators have fired at $t=t_{0}$.
This implies that in every cycle, first each neuron sends a spike before any neuron receives any spike.
Such a perturbation will affect the time of the next firing events
because the larger the perturbed phase of an oscillator is, the earlier
this oscillator reaches threshold and sends a pulse. In principle,
there one may consider other perturbations, in which, for instance,
a pulse is added or removed at a certain time. As such perturbations
are not relevant for questions of asymptotic stability, we do not
consider them here.

To construct the stroboscopic period-$T$ map, it is convenient to
rank order the elements $\delta_{i}$ of $\boldsymbol{\delta}$ in
the following manner: For each oscillator $i$ we label the perturbations
$\delta_{j}$ of its presynaptic oscillators $j\in\Pre(i)$ (for which
$\eps_{ij}\neq0$) according to their sizes\begin{equation}
\Delta_{i,1}\geq\Delta_{i,2}\geq\ldots\geq\Delta_{i,k_{i}}\label{eq:spikeorder}\end{equation}
where \begin{equation}
k_{i}:=|\Pre(i)|\end{equation}
 is the number of its presynaptic oscillators, called in-degree in
graph theory. The index $n\in\{1,\ldots,k_{i}\}$ counts the received
pulses in the order of their successive arrival. Thus, if $j_{n}\equiv j_{n}(i)\in\Pre(i)$
labels the presynaptic oscillator from which $i$ receives its $n^{\mathrm{th}}$
signal during the period considered, we have \begin{equation}
\Delta_{i,n}=\delta_{j_{n}(i)}\,.\label{eq:Deltaindeltajn}\end{equation}
For later convenience, we also define \begin{equation}
\Delta_{i,0}=\delta_{i}\,.\label{eq:Deltai0}\end{equation}

\begin{figure*}[!tph]
\begin{centering}\includegraphics[width=12cm]{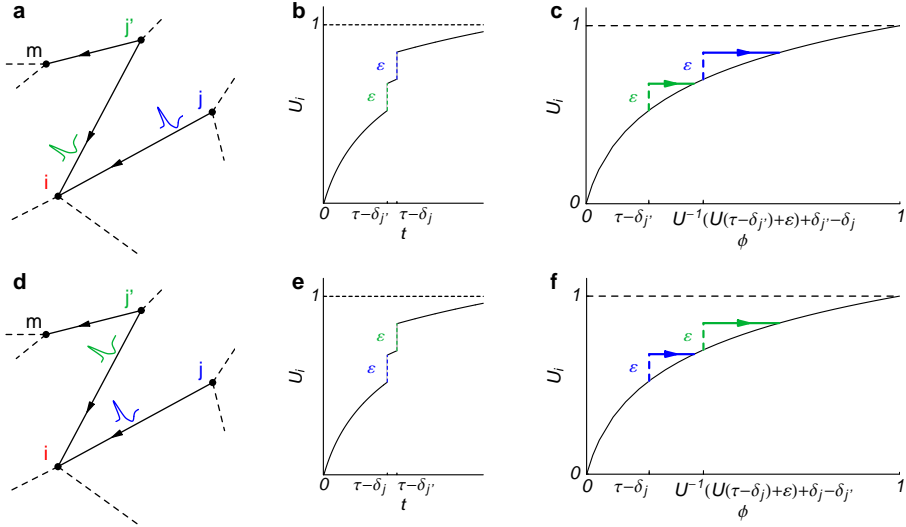}\par\end{centering}

\caption{Two signals arriving almost simultaneously induce different phase
changes, depending on their rank order. The figure illustrates a simple
case where $\Pre(i)=\{ j,j'\}$ and $\delta_{i}=0$, (a)--(c) for
$\delta_{j'}>\delta_{j}$ and (d)--(f) for $\delta_{j}>\delta_{j'\,}$.
(a), (d) Local patch of the network displaying the reception times
of signals that are received by oscillator $i$. Whereas in (a) the
signal from $j'$ arrives before the signal of $j$, the situation
in (d) is reversed. (b), (e) Identical coupling strengths induce identical
jumps of the \textit{potential} $U$ but (c),(f) the \textit{phase}
jumps these signals induce are different and depend on the order of
the incoming signals. For small $|\delta_{i}|\ll1$, individual phase
jumps are encoded by the $p_{i,n}$. \label{fig:flip_spikes}}
\end{figure*}

For illustration, let us consider an oscillator $i$ that has exactly
two presynaptic oscillators $j$ and $j'$ such that $\Pre(i)=\{ j,j'\}$
and $k_{i}=2$ (Fig.~\ref{fig:flip_spikes}a,d).%
 For certain perturbations, oscillator $i$ first receives a pulse
from oscillator $j'$ and later from oscillator $j$. This determines
the rank order, $\delta_{j'}>\delta_{j}$, and hence $\Delta_{i,1}=\delta_{j'}$
and $\Delta_{i,2}=\delta_{j}$ (Fig.~\ref{fig:flip_spikes}a). Perturbations
with the opposite rank order, $\delta_{j}>\delta_{j'}$, lead to the
opposite labeling, $\Delta_{i,1}=\delta_{j}$ and $\Delta_{i,2}=\delta_{j'}$
(Fig.~\ref{fig:flip_spikes}d). 

We now consider a fixed, arbitrary, perturbation, the rank order of
which determines the $\Delta_{i,n}$ according to the inequalities
(\ref{eq:spikeorder}). Using the phase shift function $H_{\eps}(\phi)$
(see (\ref{eq:H})) and denoting \begin{equation}
D_{i,n}:=\Delta_{i,n-1}-\Delta_{i,n}\label{eq:D_i,n}\end{equation}
 for $n\in\{1,\ldots,k_{i}\}$ we calculate the time evolution of
phase-perturbations $\delta_{i}$ satisfying (\ref{eq:taubounds}).
Without loss of generality, we choose an initial condition near $\phi_{0}(0)=\tau/2$
. The stroboscopic time-$T$ map of the perturbations, $\delta_{i}\mapsto\delta_{i}(T)$,
is obtained from the scheme given in Table \ref{Tab:complexmap}.%
\begin{table}[!tph]
\begin{centering}\begin{tabular}{|c|c|}
$t$&
$\phi_{i}(t)$\tabularnewline
\hline 
$\begin{array}{c}
0\\
\tauh-\Delta_{i,1}\\
\tauh-\Delta_{i,2}\\
\vdots\\
\tauh-\Delta_{i,k_{i}}\\
\tauh-\Delta_{i,k_{i}}+1-\beta_{i,k_{i}}\end{array}$&
$\begin{array}{c}
\frac{\tau}{2}+\delta_{i}=:\frac{\tau}{2}+\Delta_{i,0}\\
\hD{\tau}{1}\\
\hD{\beta_{i,1}}{2}\\
\vdots\\
\hD{\beta_{i,k_{i}-1}}{k_{i}}\\
\mbox{reset:}\,\,1\mapsto0\end{array}$\tabularnewline
\end{tabular}\par\end{centering}

\caption{Time evolution of oscillator $i$ in response to $k_{i}$ successively
incoming signals from its presynaptic oscillators $j_{n}$, $n\in\{1,\ldots,k_{i}\}$,
from which $i$ receives the $n^{\mathrm{th}}$ signal during this
period. The right column gives the phases $\phi_{i}(t)$ at times
$t$ given in the left column. The time evolution is shown for a part
of one period ranging from $\phi_{i}\approx\tau/2$ to reset, $1\reset0$,
such that $\phi_{i}=0$ in the last row. The first row gives the initial
condition $\phi_{i}(0)=\tau/2+\delta_{i}\,$. The following rows describe
the reception of the $k_{i}$ signals during this period whereby the
phases are mapped to $\beta_{i,n}$ after the $n^{\mathrm{th}}$ signal
has been received. The last row describes the reset at threshold such
that the respective time $T_{i}^{(0)}=\tau/2-\Delta_{i,k_{i}}+1-\beta_{i,k_{i}}$
gives the time to threshold of oscillator $i$.}

\label{Tab:complexmap}
\end{table}
 The time to threshold of oscillator $i$, which is given in the lower
left of the scheme, \begin{equation}
T_{i}^{(0)}:=\frac{\tau}{2}-\Delta_{i,k_{i}}+1-\beta_{i,k_{i}}\label{eq:time_to_spike}\end{equation}
 is about $\phi_{0}(0)=\tau/2$ smaller than the period $T$. Hence
the period-$T$ map of the perturbation can be expressed as \begin{equation}
\delta_{i}(T)=T-T_{i}^{(0)}-\frac{\tau}{2}=\beta_{i,k_{i}}-\alpha+\Delta_{i,k_{i}}\label{eq:delta_i}\end{equation}
where $\alpha$ is given by Equation (\ref{eq:alpha}).

Equation (\ref{eq:delta_i}) defines a map valid for one particular
rank order of perturbations. In general, the perturbations of all
$k_{i}$ presynaptic oscillators of oscillator $i$ lead to $k_{i}!$
different possibilities of ordering. Thus the number of possible maps,
$\mu$, is bounded by\begin{equation}
\left(\max_{i}k_{i}\right)!\leq\mu\leq(N-1)!\,.\end{equation}
 Here the minimum is assumed, $\mu=\left(\max_{i}k_{i}\right)!$,
if only one oscillator has exactly $\mu$ presynaptic oscillators
and all other oscillators have exactly one presynaptic oscillator.
The maximum, $\mu=(N-1)!\,$, is assumed if the oscillators are coupled
all-to-all such that all connections are present, $\eps_{ij}\neq0$
for all $i$ and $j\neq i$. If the coupling is all-to-all and in
addition homogeneous, $\eps_{ij}=\eps/(N-1)$ for all $i$ and $j\neq i$,
all maps are equivalent in the sense that for any pair of maps there
is a permutation of oscillator indices that transforms one map onto
the other. For general network connectivities, however, there is no
such permutation equivalence. For instance, in random networks with
$N$ vertices and edges that are independently chosen with identical
probability $p$, the number of maps increases strongly with $N$.

\section{Multi-Operator Dynamics of Small Perturbations}

\label{sec:First-order-operators}

In order to perform a local stability analysis, we consider the first
order approximations of the maps derived in the previous section.
Expanding $\beta_{i,k_{i}}$ for small $D_{i,n}\ll1$ one can proof
by induction (see Appendix) that to first order \begin{equation}
\beta_{i,k_{i}}\doteq\alpha+\sum_{n=1}^{k_{i}}p_{i,n-1}D_{i,n}\label{eq:beta_i}\end{equation}
 where \begin{equation}
p_{i,n}:=\frac{U'(U^{-1}(U(\tau)+\sum_{m=1}^{n}\eps_{ij_{m}}))}{U'(U^{-1}(U(\tau)+\eps))}\label{eq:fractions}\end{equation}
 for $n\in\{0,\,1,\,\ldots,k_{i}\}$ encodes the effect of an individual
incoming signal of strength $\eps_{ij_{n}}$. The statement $x\doteq y$
means that $x=y+\sum_{i,n}\mathcal{O}(D_{i,n}^{2})$ as all $D_{i,n}\rightarrow0$.
Substituting the first order approximation (\ref{eq:beta_i}) into
(\ref{eq:delta_i}) using (\ref{eq:D_i,n}) leads to\begin{eqnarray}
\delta_{i}(T)\doteq\sum_{n=1}^{k_{i}}p_{i,n-1}(\Delta_{i,n-1}-\Delta_{i,n})+\Delta_{i,k_{i}}\label{eq:delta_i(T)_for_sum}\end{eqnarray}
such that after rewriting \begin{equation}
\delta_{i}(T)\doteq p_{i,0}\Delta_{i,0}+\sum_{n=1}^{k_{i}}(p_{i,n}-p_{i,n-1})\Delta_{i,n}\label{eq:delta_i(T)_for_sum2}\end{equation}
to first order in all $\Delta_{i,n}$. Since $\Delta_{i,n}=\delta_{j_{n}(i)}$
for $n\in\{1,\ldots,k_{i}\}$ and $\Delta_{i,0}=\delta_{i}$ according
to Eqs. (\ref{eq:Deltaindeltajn}) and (\ref{eq:Deltai0}), this represents
a linear map\begin{equation}
\boldsymbol{\delta}(T)\doteq A\boldsymbol{\delta}\label{eq:matrixequation}\end{equation}
 where the elements of the matrix $A$ are given by \begin{equation}
A_{ij}=\left\{ \begin{array}{ll}
p_{i,n}-p_{i,n-1} & \mbox{if}\, j=j_{n}\in\Pre(i)\\
p_{i,0} & \mbox{if}\, j=i\\
0 & \mbox{if}\, j\notin\Pre(i)\cup\{ i\}.\end{array}\right.\label{eq:matrixelements}\end{equation}
 Because $j_{n}$ in (\ref{eq:matrixelements}) identifies the $n^{\mathrm{th}}$
pulse received during this period by oscillator $i$, the first order
operator depends on the rank order of the perturbations, $A=A(\rank(\boldsymbol{\delta}))$.
The map $A\boldsymbol{\delta}$ consists of a number of linear operators,
the domains of which depend on the rank order of the specific perturbation.
Thus $A\boldsymbol{\delta}$ is piecewise linear in $\boldsymbol{\delta}$.
This map is continuous but not in general differentiable at the domain
boundaries where $\delta_{i}=\delta_{j}$ for at least one pair $i$
and $j$ of oscillators. In general, signals received at similar times
by the same oscillator induce different phase changes: For the above
example of an oscillator $i$ with exactly two presynaptic oscillators
$j$ and $j'$ and equal coupling strengths, $\eps_{i,j}=\eps_{i,j'}\,$,
the first of the two received signals has a larger effect than the
second, by virtue of the concavity of $U(\phi)$. For small $|\delta_{i}|\ll1$,
this effect is encoded by the $p_{i,n}$ (see Eq.~(\ref{eq:fractions})
and the Appendix for details). Since the matrix elements Eq. (\ref{eq:matrixelements})
are differences of these $p_{i,n}$ the respective matrix elements
$A_{i,j}$ and $A_{i,j'}$ have in general different values depending
on which signal is received first. This is induced by the structure
of the network in conjunction with the pulse coupling. For networks
with homogeneous, global coupling different matrices $A$ can be identified
by an appropriate permutation of the oscillator indices. In general,
however, this is impossible.

\section{Bounds on the eigenvalues}

\label{Sec:Bounds_on_Eigenvalues}

The above consideration establish that for many network structures
when considering the stability of the synchronous state one is faced
with the task of characterizing a large number of operators instead
of a single stability matrix. Fortunately it is possible to characterize
spectral properties common to all operators by studying bounds on
their eigenvalues.

\subsection{General properties of matrix elements}

It is important to observe that the matrix elements defined by Eqs.~(\ref{eq:matrixequation})
and (\ref{eq:matrixelements}) are normalized row-wise,\begin{eqnarray}
\sum_{j=1}^{N}A_{ij} & = & A_{ii}+\sum_{\mysubstack{j=1}{j\neq i}}^{N}A_{ij}\label{eq:sum_A_ij_1}\\
 & = & A_{ii}+\sum_{n=1}^{k_{i}}A_{ij_{n}}\label{eq:fedcba}\\
 & = & p_{i,0}+\sum_{n=1}^{k_{i}}(p_{i,n}-p_{i,n-1})\\
 & = & p_{i,k_{i}}\\
 & = & 1\end{eqnarray}
 for all $i$. Here the second last equality holds because the telescope
sum equals $p_{i,k_{i}}-p_{i,0}\,$. Therefore, every matrix $A$
has the trivial eigenvalue \begin{equation}
\lambda_{1}=1\label{eq:lambda11}\end{equation}
 corresponding to the eigenvector \begin{equation}
v_{1}=(1,1,\ldots,1)^{\mathsf{T}}\label{eq:v1_1111}\end{equation}
reflecting the time-translation invariance of the system. In addition,
the diagonal elements \begin{equation}
A_{ii}=p_{i,0}=\frac{U'(\tau)}{U'(U^{-1}(U(\tau)+\eps))}=:A_{0}\label{eq:A_ii}\end{equation}
are identical for all $i$. Since $U$ is monotonically increasing,
$U'(\phi)>0$ for all $\phi$, the diagonal elements are positive,
\begin{equation}
A_{0}>0\,.\label{eq:A0_positive}\end{equation}
One should note that the matrices $A$ have the properties (\ref{eq:sum_A_ij_1})--(\ref{eq:A0_positive})
independent of single neuron parameters, the network connectivity,
and the specific perturbation considered. Due to the above properties
of the stability matrices, bounds on the eigenvalues of a specific
matrix $A$ can be obtained from the Gershgorin theorem \cite{Gershgorin:1931:749}
(see also \cite{Stoer:1992}). 

\begin{mytheorem}[Gershgorin] Given an $N\times N$ matrix $A=(A_{ij})$
and disks\begin{equation}
K_{i}:=\{ z\in\mathbb{C}|\,\,|z-A_{ii}|\leq\sum_{\mysubstack{j=1}{j\neq i}}^{N}|A_{ij}|\}\end{equation}
for $i\in\{1,\ldots,N\}$. Then the union \begin{equation}
K:=\bigcup_{i=1}^{N}K_{i}\end{equation}
 contains all eigenvalues of $A$.

\end{mytheorem}

Let us remark that for real matrices $A$ the disks $K_{i}$ in the
complex plane are centered on the real axis at $A_{ii}=A_{0}$.

\subsection{Eigenvalues for inhibitory coupling}

\begin{figure}[!tph]
\begin{centering}\includegraphics[width=8cm]{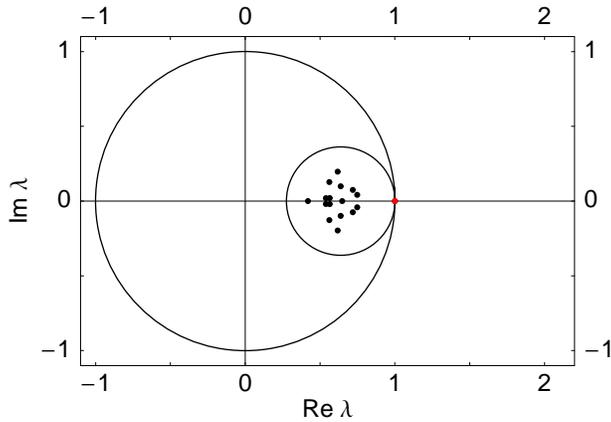}\par\end{centering}

\caption{Eigenvalues in the complex plane for inhibitory coupling. The Gershgorin
disk, that contacts the unit disk from the inside, contains all eigenvalues
of the stability matrices $A$. Black dots show eigenvalues of a specific
stability matrix for a network of $N=16$ oscillators (in which every
connection is present with probability $p=0.25)$ with coupling strengths
$\eps_{ij}=\eps/k_{i}\,$, $\eps=-0.2$, $\tau=0.15$, and a particular
rank order of the perturbation.\label{Fig:EWsample:inh}}
\end{figure}

As an application of this theorem to the above eigenvalue problem,
we consider the class of networks of inhibitorily coupled oscillators,
where all $\eps_{ij}\leq0$ and $\eps<0$. In these cases, all nonzero
matrix elements $A_{ij}$ are positive: Since $U(\phi)$ is monotonically
increasing, $U'>0$, and concave (down), $U''<0$, its derivative
$U'$ is positive and monotonically decreasing. Thus all $p_{i,n}$
(Eq.~(\ref{eq:fractions})) are positive, bounded above by one,\begin{equation}
0<p_{i,n}\leq1,\end{equation}
and increase with $n$,\begin{equation}
p_{i,n-1}<p_{i,n}\,.\end{equation}
 Hence the nonzero off-diagonal elements are positive, $A_{ij_{n}}=p_{i,n}-p_{i,n-1}>0$
such that\begin{equation}
A_{ij}\geq0\label{eq:A_ij_inh_positive}\end{equation}
for all $i,j\in\{1,\ldots,N\}$ and \begin{equation}
0<A_{0}<1.\label{eq:A0inh}\end{equation}
As a consequence, for inhibitorily coupled oscillators, \begin{equation}
\sum_{\mysubstack{j=1}{j\neq i}}^{N}|A_{ij}|=\sum_{j=1}^{N}A_{ij}-A_{ii}=1-A_{0}\end{equation}
 such that all Gershgorin disks $K_{i}$ are identical and the disk
\begin{equation}
K_{i}=K=\{ z\in\mathbb{C}|\,\,|z-A_{0}|\leq1-A_{0}\}\label{eq:Kinh}\end{equation}
 contains all eigenvalues of $A$. This disk $K$ contacts the unit
disk from the inside at the trivial eigenvalue $z=\lambda_{1}=1$
(Fig.~\ref{Fig:EWsample:inh}).%

Since the unit circle separates stable from unstable eigenvalues and
structural perturbations may move any but the trivial unit eigenvalue,
it is interesting to examine whether the unite eigenvalues may be
degenerate. For strongly connected networks the Perron-Frobenius theorem
implies that this eigenvalue is unique demonstrating the structural
stability of the confinement of eigenvalues to the unit circle. 

If the network is strongly connected, the resulting stability matrix
$A$ is irreducible such that the Perron-Frobenius theorem \cite{Perron:1907:248,Frobenius:1908:471,Frobenius:1909:514,Frobenius:1912:456}
(see also \cite{Minc:1988,Mehta:1989}) is applicable (we only state
the theorem partially).

\begin{mytheorem}[Perron-Frobenius] Let $A$ be an $N\times N$
irreducible matrix with all its elements real and non-negative. Then

\begin{itemize}
\item $A$ has a real positive eigenvalue $\lambda_{\textrm{max}}$, the
maximal eigenvalue, which is simple and such that all eigenvalues
$\lambda_{i}$ of $A$ satisfy $|\lambda_{i}|\leq\lambda_{\textrm{max}}$.
\end{itemize}
\end{mytheorem}

The Perron-Frobenius theorem implies, that the eigenvalue that is
largest in absolute value, here $\lambda_{1}=1$, is unique for strongly
connected networks. The Gershgorin theorem guarantees that eigenvalues
$\lambda$ of modulus one are degenerate, $\lambda=\lambda_{1}=1$.
Taken together, for strongly connected networks, the non-trivial eigenvalues
$\lambda_{i}$ satisfy \begin{equation}
|\lambda_{i}|<1\label{eq:lambdasmallerone}\end{equation}
 for $i\in\{2,\ldots,N\}$. This suggests that the synchronous state
is stable for inhibitory couplings. As pointed out below ( \sect{stability}
a proof of stability, however, requires further analysis.

\subsection{Eigenvalues for excitatory coupling}

\begin{figure}[!tph]
\begin{centering}\includegraphics[width=8cm]{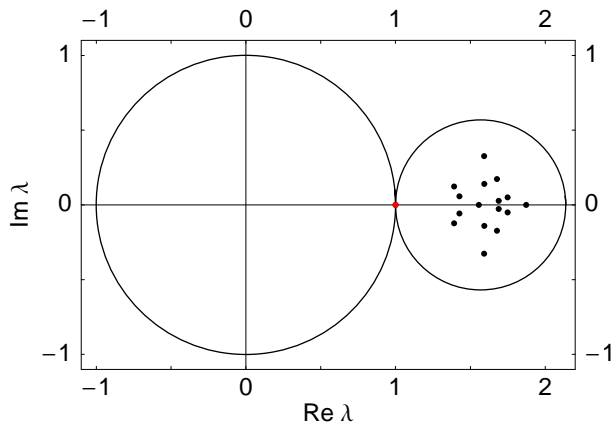}\par\end{centering}

\caption{Eigenvalues in the complex plane for excitatory coupling. Gershgorin
disk, that contacts the unit disk from the outside, contains all eigenvalues
of the stability matrices $A$. Black dots show eigenvalues of a specific
stability matrix for a network of $N=16$ oscillators (in which every
connection is present with probability $p=0.25)$ with coupling strength
$\eps_{ij}=\eps/k_{i}\,$,$\eps=0.2$, $\tau=0.15$) and a particular
rank order of the perturbation.\label{Fig:EWsample:exc}}
\end{figure}

Let us briefly discuss the case of excitatorily coupled oscillators,
where all $\eps_{ij}\geq0$ and $\eps>0$. Here the analysis proceeds
similar to the case of inhibitory coupling. Due to the monotonicity
and concavity of $U(\phi)$, we obtain\begin{equation}
p_{i,n}\geq1\end{equation}
as well as a decrease with $n$,\begin{equation}
p_{i,n-1}>p_{i,n}\,,\end{equation}
such that $A_{i,j_{n}}=p_{i,n}-p_{i,n-1}<0$ and thus\begin{equation}
A_{ij}\leq0\end{equation}
for $i\in\{1,\ldots,N\}$ and $j\neq i$ and \begin{equation}
A_{ii}=A_{0}>1.\end{equation}
 Consequently, for excitatorily coupled oscillators,\begin{equation}
\sum_{\mysubstack{j=1}{j\neq i}}^{N}|A_{ij}|=-\sum_{\mysubstack{j=1}{j\neq i}}^{N}A_{ij}=A_{0}-1\end{equation}
such that again all Gershgorin disks $K_{i}$ are identical and \begin{equation}
K_{i}=K=\{ z\in\mathbb{C}|\,\,|z-A_{0}|\leq A_{0}-1\}.\end{equation}
Since $A_{0}>1$, the disk $K$ contacts the unit disk from the outside
at $z=\lambda_{1}=1$ (Fig.~\ref{Fig:EWsample:exc}).%
 If the network is strongly connected, $\lambda_{1}=1$ is again unique
by the Perron-Frobenius theorem for irreducible matrices because it
is the largest eigenvalue (in absolute value) of the inverse $A^{-1}$
of $A$, if the inverse exists. This result indicates that the fully
synchronous state is unstable for excitatory couplings.

\section{Stability\Lsect{stability}}

\label{sec:Stability}

In Sec.\ \ref{Sec:Bounds_on_Eigenvalues}, we found analytical bounds
on the eigenvalues of the stability matrices. However, even if only
eigenvalues $\lambda$ with $|\lambda|\leq1$ are present, a growth
of perturbations might seem to be possible because (i) the eigenspaces
of the (asymmetric) matrices $A$ cannot be guaranteed to be orthogonal,
and (ii) in general, different matrices $A$ are successively applied
to a given perturbation. Due to the non-orthogonality (i) the length
$\left\Vert \boldsymbol{\delta}\right\Vert $ of a given perturbation
vector $\boldsymbol{\delta}$ might increase during one period. Since
the stability matrix and the set of eigenvectors may change due to
(ii), the length of the perturbation vector might increase in the
subsequent period as well. Since this procedure may be iterated, the
eigenvalues $\lambda_{i}$ of the stability matrices, although satisfying
$|\lambda_{i}|<1$ for $i\in\{2,\ldots,N\}$ and $\lambda_{1}=1$,
are not sufficient to ensure the stability of the synchronous state. 

Thus, the eigenvalues of the dynamically changing stability matrices
guide the intuition about the stability of the synchronous state as
well as about the speed of convergence (in case of stability) or divergence
(in case of instability). Nevertheless, stability cannot be directly
inferred from the set of eigenvalues. In the following we illustrate
the final proof of stability in the simple case of inhibitory coupling
where all $A_{ij}\geq0$ (\ref{eq:A_ij_inh_positive}).

\subsection{Plain stability for inhibitory coupling}

For inhibitory networks, the proof of plain (non-asymptotic) stability
is simple. Given the fact that for inhibition $\sum_{j=1}^{N}A_{ij}=1$
and $A_{ij}\geq0$, the synchronous state is stable because a given
perturbation $\boldsymbol{\delta}=\boldsymbol{\delta}(0)$ satisfies

\begin{eqnarray}
\left\Vert \boldsymbol{\delta}(T)\right\Vert  & := & \max_{i}|\delta_{i}(T)|\\
 & = & \max_{i}\left|\sum_{j=1}^{N}A_{ij}\delta_{j}\right|\\
 & \leq & \max_{i}\sum_{j}|A_{ij}||\delta_{j}|\\
 & \leq & \max_{i}\sum_{j}|A_{ij}|\max_{k}|\delta_{k}|\\
 & = & \max_{i}\sum_{j}A_{ij}\max_{k}|\delta_{k}|\\
 & = & \max_{k}|\delta_{k}|\\
 & = & \left\Vert \boldsymbol{\delta}\right\Vert \,.\end{eqnarray}
In this section, we use the vector norm \begin{equation}
\left\Vert \boldsymbol{\delta}\right\Vert :=\max_{i}\delta_{i}\,.\end{equation}
 Thus the length of a perturbation vector can not increase during
one period implying that it does not increase asymptotically. This
result is independent of the connectivity structure of the network,
the special choice of parameters, $\eps_{ij}\leq0$, $\tau>0$, the
potential function $U(\phi)$, and the rank order of the perturbation.

\subsection{Asymptotic stability in strongly connected networks}

Asymptotic stability can be established for networks of inhibitorily
coupled oscillators assuming that the network satisfies the condition
of strong connectivity. A directed graph is called strongly connected
if that every vertex can be reached from every other vertex by following
the directed connections. Thus a network is strongly connected, if
every oscillator can communicate with each oscillator in the network
at least indirectly. It is clear that in a disconnected network only
the connected components may synchronize completely in the long-term,
but these components can not be synchronized by mutual interactions.
In the proof given below, we do not consider networks that are disconnected.
We do also not consider networks that are weakly connected (but not
strongly connected) such as two globally coupled sub-networks which
are linked by unidirectional connections from one sub-network to the
other. 

The synchronous state may be characterized by a perturbation $\boldsymbol{\delta}\equiv\boldsymbol{\delta}(0)$
that represents a uniform phase shift, \begin{equation}
\boldsymbol{\delta}=c_{1}\mathbf{v}_{1}=c_{1}(1,1,\ldots,1)^{\textsf{T}}\label{eq:c1v1}\end{equation}
 where $c_{1}\in\mathbb{R}$, $|c_{1}|\ll1$, and $\mathbf{v}_{1}$
is the eigenvector of $A$ corresponding to the eigenvalue $\lambda_{1}=1$
(\ref{eq:v1_1111}). Such a perturbation satisfies \begin{equation}
\boldsymbol{\delta}(T)=A\boldsymbol{\delta}=\boldsymbol{\delta}\label{eq:deltav1invariant}\end{equation}
for all matrices $A=A(\rank(\boldsymbol{\delta}))$ independent of
the rank order $\rank(\boldsymbol{\delta})$ and thus\begin{equation}
\left\Vert \boldsymbol{\delta}(T)\right\Vert =\left\Vert \boldsymbol{\delta}\right\Vert \,.\end{equation}
Now consider a $\boldsymbol{\delta}$ that does not represent the
synchronous state, \begin{equation}
\boldsymbol{\delta}\neq c_{1}\mathbf{v}_{1}\label{eq:delta_not_v1}\end{equation}
for all $c_{1}\in\mathbb{R}$. Then one might guess that\begin{equation}
\left\Vert \boldsymbol{\delta}(T)\right\Vert <\left\Vert \boldsymbol{\delta}\right\Vert .\end{equation}
 We assume that (i) $\boldsymbol{\delta}$ is not the synchronous
state (\ref{eq:delta_not_v1}) and that (ii) all single-oscillator
perturbations are non-negative, $\delta_{i}\geq0$. The latter assumption
is made without loss of generality, because otherwise a vector proportional
to $\mathbf{v}_{1}$ can be added such that $\delta_{i}\geq0$ is
satisfied. The assumption (ii) implies that the components of the
perturbation stay non-negative for all times, because $\delta_{i}(T)=\sum_{j=1}^{N}A_{ij}\delta_{j}\geq0$
for all $i$ such that $\delta_{i}(lT)\geq0$ for all $i$ and all
$l\in\NN$. 

For convenience, we define the largest component of the perturbation\begin{equation}
\delta_{M}:=\max_{i}\delta_{i}\end{equation}
and the second largest component \begin{equation}
\delta_{m}:=\max\{\delta_{i}\,|\,\,\delta_{i}<\delta_{M}\}\end{equation}
such that $\delta_{m}<\delta_{M}$ . We also define the index set
of maximal components \begin{equation}
\mathbb{M}:=\{ j\in\{1,\ldots,N\}\,|\,\delta_{j}=\delta_{M}\}\end{equation}
which is always non-empty. We write $j\notin\mathbb{M}$ if $j\in\{1,\ldots,N\}\setminus\mathbb{M}$.
With these definitions we find\begin{eqnarray}
\delta_{i}(T) & = & \sum_{j=1}^{N}A_{ij}\delta_{j}\\
 & = & \sum_{j\in\mathbb{M}}A_{ij}\delta_{j}+\sum_{j\notin\mathbb{M}}A_{ij}\delta_{j}\\
 & \leq & \sum_{j\in\mathbb{M}}A_{ij}\delta_{M}+\sum_{j\notin\mathbb{M}}A_{ij}\delta_{m}\\
 &  & +\sum_{j\notin\mathbb{M}}A_{ij}\delta_{M}-\sum_{j\notin\mathbb{M}}A_{ij}\delta_{M}\\
 & = & \delta_{M}\sum_{j=1}^{N}A_{ij}-(\delta_{M}-\delta_{m})\sum_{j\notin\mathbb{M}}A_{ij}\\
 & = & \delta_{M}-(\delta_{M}-\delta_{m})\sum_{j\notin\mathbb{M}}A_{ij}\,.\label{eq:asymptoticprooflastrow}\end{eqnarray}

Hence if $\sum_{j\notin\mathbb{M}}A_{ij}>0$ the norm of the perturbation
vector \textit{decreases} in one period,\begin{equation}
\left\Vert \boldsymbol{\delta}(T)\right\Vert =\max_{i}\delta_{i}(T)<\delta_{M}=\max_{i}\delta_{i}=\left\Vert \boldsymbol{\delta}\right\Vert .\label{eq:delta_max(T)_lt_delta_max(0)}\end{equation}

There are, however, also perturbations that imply $A_{ij}=0$ for
(at least) one specific $i$ and all $j\notin\mathbb{M}$ such that
$\sum_{j\notin\mathbb{M}}A_{ij}=0$. This is the case if and only
if there is an oscillator $i$ that receives input only from oscillators
$j$ with maximal components, $\delta_{j}=\delta_{M}$, and itself
has maximal component, $\delta_{i}=\delta_{M}$ . So suppose that
\begin{equation}
\exists\, i\in\{1,\ldots,N\}\,\,\forall\, j\in\Pre(i)\cup\{ i\}\,:\,\,\delta_{j}=\delta_{M}\,,\label{eq:i_Pre_i_max}\end{equation}
i.e.\ $\Pre(i)\cup\{ i\}\subset\mathbb{M}$. Then $\delta_{i}(T)=\delta_{M}$
for this oscillator $i$ and hence \begin{equation}
\max_{i}\delta_{i}(T)=\delta_{M}\end{equation}
such that the norm of the perturbation vector does \textit{not decrease}
within one period.

Nevertheless, if the network is strongly connected, the norm of the
perturbation vector is reduced in at most $N-1$ periods: We know
that if and only if there is no oscillator $i$ satisfying (\ref{eq:i_Pre_i_max}),
the norm will be reduced (\ref{eq:delta_max(T)_lt_delta_max(0)})
because $\sum_{j\notin\mathbb{M}}A_{ij}>0$ for all $i$ by (\ref{eq:asymptoticprooflastrow}).
So to have $\max_{i}\delta_{i}(T)=\delta_{M}$ one needs one oscillator
$i$ that satisfies (\ref{eq:i_Pre_i_max}), i.e.\ $i$ and $\Pre(i)$
have to have maximal components. If the vector norm stays maximal
another period,\begin{equation}
\max_{i}\delta_{i}(2T)=\delta_{M}\,,\end{equation}
 not only all $j\in\{ i\}\cup\Pre(i)$ but also all $j\in\Pre(\Pre(i))$
have to have maximal components, $\delta_{j}=\delta_{M}$. Iterating
this $l$ times, leads to the condition \begin{eqnarray}
\begin{array}{c}
\max_{i}\delta_{i}(lT)=\delta_{M}\,\Leftrightarrow\\
\\\exists\, i\,\,\forall\, j\in\{ i\}\cup\Pre(i)\cup\Pre^{(2)}(i)\cup\ldots\cup\Pre^{(l)}(i)\,:\,\,\delta_{j}=\delta_{M}\end{array}\end{eqnarray}
where \begin{equation}
\Pre^{(l)}(i):=\underbrace{{\Pre\circ\Pre\circ\ldots\circ\Pre}}_{l\textrm{ times}}(i)\end{equation}
 is the set of oscillators, that is connected to oscillator $i$ via
a sequence of exactly $l$ directed connections.

Since for a strongly connected network the union of all presynaptic
oscillators and their respective presynaptic oscillators is the set
of all oscillators \begin{equation}
\{ i\}\cup\Pre(i)\cup\Pre^{(2)}(i)\cup\ldots\cup\Pre^{(l)}(i)=\{1,\ldots,N\}\label{eq:all_presynpses_all}\end{equation}
for $l\geq N-1$, this leads to the conclusion that \begin{equation}
\max_{i}\delta_{i}((N-1)T)=\delta_{M}\,\Rightarrow\,\forall j\,:\,\delta_{j}=\delta_{M}\end{equation}
 such that \begin{equation}
\boldsymbol{\delta}=\delta_{M}\mathbf{v}_{1}\end{equation}
 contrary to the assumption (\ref{eq:delta_not_v1}). Note that for
a given network connectivity, the condition (\ref{eq:all_presynpses_all})
is satisfied for any $l\geq l_{c}$ where $l_{c}$ is the diameter
of the underlying network, the longest directed connection path between
any two oscillators in the network. The diameter is maximal, $l_{c}=N-1$,
assumed for a ring of $N$ oscillators.

Hence, after at most $l_{c}$ periods, the norm of the perturbation
vector decreases,\begin{equation}
\left\Vert \boldsymbol{\delta}(lT)\right\Vert <\left\Vert \boldsymbol{\delta}\right\Vert \end{equation}
for $l\geq l_{c}$. Since $\left\Vert \boldsymbol{\delta}(ml_{c}T)\right\Vert $
is strictly monotonically decreasing with $m\in\mathbb{N}$ and bounded
below by zero, the limit \begin{equation}
\boldsymbol{\delta}^{\infty}:=\lim_{m\rightarrow\infty}\boldsymbol{\delta}(ml_{c}T)\end{equation}
exists. If $\boldsymbol{\delta}^{\infty}\neq c_{1}\mathbf{v}_{1}$
the vector norm would be reduced as implied by the above considerations.
Thus we find that\begin{equation}
\boldsymbol{\delta}^{\infty}=c_{1}\mathbf{v}_{1}\end{equation}
because the uniform components $c_{1}\mathbf{v}_{1}$ of the original
perturbation $\boldsymbol{\delta}$ does not change under the dynamics
(see Eq.~(\ref{eq:deltav1invariant})). This completes the demonstration
that the synchronous state is asymptotically stable for any strongly
connected network of inhibitorily coupled oscillatory units.

\section{Summary and discussion}

\label{sec:Conclusions}

We analyzed the stability of synchronous states in arbitrarily connected
networks of pulse-coupled oscillators. For generally structured networks,
the intricate problem of multiple stability operators arises. We analyzed
this multi-operator problem exactly. For both inhibitory and excitatory
couplings, we determined analytical bounds on the eigenvalues of all
operators. Given the multi-operator property, stability cannot be
deduced by considering the eigenvalues only. We therefore completed
the stability analysis on two levels For networks with inhibitory
couplings ($\eps_{ij}\leq0$) we found plain (Lyapunov) stability
of the synchronous state; under mild constraints (strong connectivity
of the network), graph theoretical arguments show that it is also
asymptotically stable, independent of the parameters and the details
of the connectivity structure. 

Let us point out that the stability results obtained here are valid
for arbitrarily large finite networks of general connectivity. In
contrast to ordinary stability problems, the eigenvalues of the first
order operators do not directly determine the stability here, because
different linear operators may act subsequently on a given perturbation
vector. Thus the eigenvalues and eigenvectors of an individual matrix
alone do not define the local dynamics of the system. 
For network dynamical systems of pulse-coupled units that exhibit more complex cluster states, with two or more groups of synchronized units, an analogous multi-operator problem arises \cite{Timme:2002:154105,Timme:2003:1,Ashwin:2005:2035}. Methods to bound their spectra and the graph theoretical approach presented above to prove asymptotic stability are applicable to these cluster states in a similar way, provided the stability problem is originally formulated in an event-based manner. Further studies
(e.g., \cite{Memmesheimer:2006:188101}) show that the stable synchrony
in complex networks is not restricted to the class of models considered
here. Together with other works, e.g. \cite{Denker:2002:074103,Timme:2004:074101,Timme:2006:015108},
this also indicates robustness of the synchronous state against structural
perturbations, \ie conditions necessary to ensure that the local
dynamics stays qualitatively similar if the time evolution of the
model is weakly modified. 

It has been hypothesized (see \eg \cite{Abeles:1991,Diesmann:1999:529})
that the experimentally observed precisely timed firing patterns in
otherwise irregular neural activity might be generated dynamically
by the cerebral cortex. As a by-product, our results clearly demonstrate
that states in which units fire in a temporally highly regular fashion
and states with irregular asynchronous activity may be coexisting
attractors of the same recurrent network, cf. \cite{Timme:2002:258701}. This already applies to
random networks. More specifically structured networks may possess
a large variety of dynamical states in which firing times are precisely
coordinated. A promising direction for future research is thus to
adopt the methods developed here for investigating the stability of
such states. In particular, methods adapted from the ones developed
here may reveal dynamical features of networks in which the coupling
strengths are not only structured but highly heterogeneous and precise
temporal firing patterns occur in place of a simple synchronous state
\cite{Denker:2002:074103,Memmesheimer:2006:188101,Jahnke:2007:1}.
More generally, our methods may help to understand properties of stability
and robustness in various network dynamical systems with pulsed interactions.
They may thus be relevant not only for networks of nerve cells in
the brain coupled via chemical synapses but also for cells in heart
tissue coupled electrically and even for populations of fireflies
and chirping crickets that interact by sending and receiving light
pulses and sound signals, respectively.

Acknoledgements: We thank Theo Geisel, Michael Denker and Raoul-Martin
Memmesheimer for fruitful discussions. This work has been partially
supported by the Ministry of Education and Research (BMBF) Germany
under Grant No. 01GQ07113 and, via the Bernstein Center for Computational
Neuroscience (BCCN) G\"{o}ttingen, under Grant No. 01GQ0430.

\section{Appendix: Exact derivation of the expansion (\ref{eq:beta_i})}

We prove a generalization of the first order expansion (\ref{eq:beta_i})
used above,\begin{equation}
\beta_{i,m}\doteq\alpha_{i,m}+\sum_{n=1}^{m}p_{i,n-1,m}\, D_{i,n}\quad\mathrm{for}\quad m\in\{1,\ldots,k_{i}\},\label{eq:A:beta_i,m}\end{equation}
 by induction over $m$ for arbitrary fixed $i\in\{1,\ldots,N\}$.
As in the main text, here the statement $x\doteq y$ means that $x=y+\sum_{i,n}\mathcal{O}(D_{i,n}^{2})$
as all $D_{i,n}\rightarrow0$. The quantities $\beta_{i,m}$ were
defined in Table \ref{Tab:complexmap} . The quantities (\ref{eq:A:beta_i,m})
appear in the derivation of the map (\ref{eq:delta_i}) for $m=k_{i}$.
As an extension of the notation in section \ref{sec:strobmap}, we
denote here \begin{equation}
\alpha_{i,m}=U^{-1}\left(U\left(\tau\right)+\sum_{n=1}^{m}\eps_{ij_{n}}\right)\quad\mathrm{for}\quad m\in\{0,\ldots,k_{i}\},\end{equation}
\begin{equation}
D_{i,n}=\Delta_{i,n-1}-\Delta_{i,n}\quad\mathrm{for}\quad n\in\{1,\ldots,k_{i}\}\,,\end{equation}
 and \begin{equation}
p_{i,n,m}=\frac{U'(\alpha_{i,n})}{U'(\alpha_{i,m})}\quad\mathrm{for}\quad0\leq n\leq m\leq k_{i}\,.\end{equation}
Here the prime denotes the derivative of the potential function $U$
with respect to its argument. The latter definition implies the identity\begin{equation}
p_{i,n,l}\, p_{i,l,m}=p_{i,n,m}\,.\label{eq:A:p_transitiv}\end{equation}
 For the case $m=1$, the induction basis, expression (\ref{eq:A:beta_i,m})
holds because\begin{eqnarray}
\beta_{i,1} & = & U^{-1}(U(\tau+D_{i,1})+\eps_{ij_{1}})\\
 & \doteq & U^{-1}(U(\tau)+\eps_{ij_{1}})\\
 &  & +\left[\frac{\partial}{\partial D}U^{-1}(U(\tau+D)+\eps_{ij_{1}})\right]_{D=0}D_{i,1}\nonumber \\
 & = & \alpha_{i,1}+\frac{U'(\tau)}{U'(U^{-1}(U(\tau)+\eps_{ij_{1}}))}D_{i,1}\\
 & = & \alpha_{i,1}+\frac{U'(\alpha_{i,0})}{U'(\alpha_{i,1})}D_{i,1}\\
 & = & \alpha_{i,1}+p_{i,0,1}\, D_{i,1}\end{eqnarray}
where the third last equality holds because the derivatives of inverse
functions are given by $\partial/\partial x\, U^{-1}(x)=1/U'(U^{-1}(x))$.
The induction step, $m\mapsto m+1$, is proven by \begin{eqnarray}
\beta_{i,m+1} & = & U^{-1}(U(\beta_{i,m}+D_{i,m+1})+\eps_{ij_{m+1}})\\
 & \doteq & U^{-1}\left(U\left(\alpha_{i,m}+\sum_{n=1}^{m}p_{i,n-1,m}D_{i,n}+D_{i,m+1}\right)+\eps_{ij_{m+1}}\right)\\
 & \doteq & U^{-1}\left(U\left(\alpha_{i,m}\right)+\eps_{ij_{m+1}}\right)\\
 &  & \quad+\left[\frac{\partial}{\partial D}U^{-1}(U(\alpha_{i,m}+D)+\eps_{ij_{m+1}})\right]_{D=0}\left(\sum_{n=1}^{m}p_{i,n-1,m}\, D_{i,n}+D_{i,m+1}\right)\nonumber \\
 & = & U^{-1}\left(U\left(\alpha_{i,m}\right)+\eps_{ij_{m+1}}\right)+\frac{U'(\alpha_{i,m})}{U'(\alpha_{i,m+1})}\left(\sum_{n=1}^{m}p_{i,n-1,m}\, D_{i,n}+D_{i,m+1}\right)\\
 & = & \alpha_{i,m+1}+p_{i,m,m+1}\left(\sum_{n=1}^{m}p_{i,n-1,m}\, D_{i,n}+D_{i,m+1}\right)\\
 & = & \alpha_{i,m+1}+\sum_{n=1}^{m+1}p_{i,n-1,m+1}\, D_{i,n}\end{eqnarray}
where we used $U(\alpha_{i,m})+\eps_{ij_{m+1}}=U(\alpha_{i,m+1})$
in the second last step and the identity (\ref{eq:A:p_transitiv})
in the last step. This completes the proof of the first order expansion
(\ref{eq:A:beta_i,m}). 

Note that because of the normalization $\sum_{n=1}^{k_{i}}\eps_{ij_{n}}=\sum_{j=1}^{N}\eps_{ij}=\eps$
for all $i$, the quantity $\alpha_{i,k_{i}}=U^{-1}(U(\tau)+\eps)=\alpha$
is independent of $i$. In addition, $p_{i,n,k_{i}}=p_{i,n}$ for
all $i$ and all $n$. Hence, theorem (\ref{eq:A:beta_i,m}) in the
case $m=k_{i}$ yields expression (\ref{eq:beta_i}).

\end{document}